 \documentclass[12pt]{article}
 \usepackage{amsmath,amsfonts,theorem}
 \numberwithin{equation}{section}

 \newtheorem{prop}{Proposition}[section]

 \theorembodyfont{\romfamily}

 \newcommand{\hidot}{^{\scriptscriptstyle{\bullet}}}
 \newcommand{\dotlessj}{j}

 \newcommand{\Span}[1]{\left< #1 \right>}
 \newcommand{\half}{\frac12}

 \newcommand{\Ad}{\operatorname{Ad}}
 \newcommand{\id}{\operatorname{id}}
 \newcommand{\pr}{\operatorname{pr}}
 \newcommand{\Tr}{\operatorname{Tr}}

 \newcommand{\BS}{\mathrm{BS}}
 
 \newcommand{\SNW}{\mathrm{SNW}}

 \newcommand{\tensor}{\otimes}
 \newcommand{\al}{\alpha}

 \newcommand{\De}{\Delta}
 \newcommand{\Ga}{\Gamma}
 
 \newcommand{\Om}{\Omega}
 
 \newcommand{\Si}{\Sigma}

\newcommand{\ov}{\overline}

 \newcommand{\C}{\mathbb C}

 \newcommand{\R}{\mathbb R}
 \newcommand{\Z}{\mathbb Z}

 \newcommand{\Hom}{\operatorname{Hom}}

 \newcommand{\End}{\operatorname{End}}
 \renewcommand{\Im}{\operatorname{Im}}

 \newcommand{\U}{\operatorname U} % unitary group
 \newcommand{\PU}{\operatorname{PU}}
 \newcommand{\SU}{\operatorname{SU}}
 \newcommand{\fsu}{\operatorname{\mathfrak{su}}}

\begin{document}

 \title{Three mathematical faces\\ of $\SU(2)$-spin networks}
 \markright{\hfill Three faces of\/ $\SU(2)$-spin networks \quad}

 \author{Andrei Tyurin}
 \date{13th Oct 2000 }
 \maketitle

Spin networks are at the core of quantum gravity \cite{1}. We have neither
the space nor the competence to give an exhaustive list of the physical
and
philosophical interpretations of this notion (for these, see for example
\cite {2} and \cite{3}). New (and old) approaches towards creating a
finite
quantum theory of general relativity would use combinatorial expressions
in
Feynman integrals, spin networks, spin foams and others combinatorial
objects.
The intention is to get out the standard ``continuous'' geometry.
Our aim is to plug the mathematical community at large into these
procedures
as participants. For this, because of the different cultural backrounds,
we
would like to change tack: to relate discrete (combinatorial) objects to
the
standard ``continuous'' geometry. From the mathematical point of view,
relations of this type give rise to identifications between the geometry
of
varieties and combinatorial objects, as exemplified by the relation
between
Lie algebras and root systems, or toric varieties and polytopes.

The general mathematical mechanism of such ``interpretations'' could be
called the {\it analytic\/} theory of non-Abelian theta functions, since
they run completely parallel to the original classical theory of theta
functions. Recall that this classical theory has three parts:

 \begin{enumerate}
 \item One constructs a family of holomorphic functions
$\{ \theta_{\al_k}^\Om \}$ on a fixed space (to be concrete, on
$(\C^*)^g$).
Each function depends on a symmetric $g\times g$ complex matrix $\Om$ with
positive imaginary part $\Im \Om>0$, and $\al_k$ is some combinatorial
data,
the so-called ``characteristics of level $k$''. We will see that this data
is equivalent to a choice of a $\U(1)$-spin network of genus $g$ (see the
final Section~4).

\item If $\Om$ is the period matrix of a marked Riemann surface $\Si$ of
genus $g$, the functions $\{ \theta_{\al_k}^\Om \}$ form a basis of
$H^0(J\Si, \Theta^k)$, the space of holomorphic sections of the $k$th
power of the polarizing line bundle $\Theta$ on the Jacobian of $\Si$. In
particular, for $k=1$ there is just one section (up to scaling), whose
zero
set is the theta divisor of $\Si$; by Riemann's theorem, this is
birational
to the symmetric power $\Si^{[g-1]}$.

\item The final part describes how these geometric objects behave under
deformations of $\Om$. Here we have transformation rules governing changes
of marking and projective flat connections under continuous deformation of
parameters. (We would like to emphasize the best reference for these
classical projective flat connections \cite{4}.)
\end{enumerate}

Realizing this program in the case of $\SU(2)$ is the subject of current
work
\cite{5} in collaboration with C. Florentino, J. Mour\~ao and
J.P. Nunes. The
general picture is a mosaic consisting of many stones, or a many piece
jigsaw
puzzle. This paper describes one stone of the mosaic, taken from the first
part of the program; namely, we describe the combinatorial data of every
non-Abelian theta function, the analog of a theta characteristic, as an
{\it
$\SU(2)$-spin network} of genus $g$, and associate a ``continuous''
geometry
with it. For this, we must show all three different mathematical faces of
$\SU(2)$-spin networks, although possibly only the third is new.

\section{First face}
A spin network is a labeled trivalent graph $\Ga$, with edges labeled by
finite dimensional irreducible representations of $\SU(2)$ and vertices
labeled by the associated intertwiners. Traditionally, in this case, the
set
$\widehat{SU(2)}$ of irreducible representations is the nonnegative
half-integers
 \begin{equation}
 \widehat{SU(2)}=\half \Z^{+}.
 \end{equation}
However, from the combinatorial point of view, we find it convenient to
multiply these numbers by 2, and call the result {\it colors}. Let
$E(\Ga)$
be the set of edges, $V(\Ga)$ the set of vertices and $F(\Ga)=\{v\in e \}$
the set of flags, where a flag is an edge with fixed end. Clearly
$F(\Ga)\subset E(\Ga) \times V(\Ga)$, and two projections
 \begin{equation}
 e \colon F(\Ga)\to E(\Ga) \quad \text{and} \quad
 v \colon F(\Ga)\to V(\Ga)
 \end{equation}
are ramified covers of degree 2 and 3 having the same ramification locus,
consisting of pairs $v\in e$ where the edge $e$ is a {\it loop}. (Recall
that a trivalent graph can contain loops, that is, an edge whose two ends
coincide.) If $L(\Ga)$ is the set of loops in $\Ga$, and $\vert\ \vert$
the
number of elements of a finite set then
 \begin{equation}
2 \cdot \vert E(\Ga) \vert-\vert L(\Ga) \vert=\vert F(\Ga) \vert=
3 \cdot \vert V(\Ga) \vert-\vert L(\Ga)\vert.
 \end{equation}
Hence
 \begin{equation}
 \vert V(\Ga) \vert=2g-2; \quad \vert E(\Ga) \vert=3g-3,
 \end{equation}
where $g>1$ is a certain integer called the {\it genus} of $\Ga$.

Thus a spin network defines a map
 \begin{equation}
 j \colon E(\Ga)\to \widehat{SU(2)}.
 \end{equation}

Recall that for a triple of representations $j_{1}, j_{2}, j_{3}$, an
{\it intertwiner} is a trivial component of the tensor product
$j_{1} \tensor j_{2} \tensor j_{3}$. Such a component exists iff the
Clebsch--Gordan conditions
 \begin{equation}
 j_{1}+j_{2}+j_{3}\in \Z \quad\text{and}\quad \vert j_{1}-j_{2} \vert \le
 j_{3} \le j_{1}+j_{2} 
 \end{equation}
hold for every ordering of edges around a vertex.

A function $j$ (1.5) defines a {\it spin network} $\Ga_j$ iff these
conditions
hold for every triple $j_{v,1}, j_{v,2}, j_{v,3}$ of representations
around
every vertex $v\in V(\Ga)$ and for every ordering of triple of edges.
In this case any intertwiner
 \begin{equation}
 i_v\in j_{v,1} \tensor j_{v,2} \tensor j_{v,3}
 \end{equation}
is defined uniquely. Thus we can omit any labeling of the vertices and
denote
a $\SU(2))$-spin network by the symbol $\Ga_{j}$.

A spin network is of {\it level} $k$ if $j_{i} \le k$ for every edge
$e_i$. There is a finite number $N^{k}(\Ga)$ of spin networks of level $k$
with graph $\Ga$, and a finite number $N^{k}_{g}$ of all spin networks of
level $k$ and genus $g$.

Penrose prescribed a number $V(\Ga_{j})$ for every spin network $\Ga_j$,
its
{\it value}. We omit the precise definition here for reasons of space: the 
``right'' definition involves starting with loop representations, as in
the
beautiful paper \cite{6}. 

A one-trivalent graph $\Ga$ (or a {\it history}) is a trivalent graph
possibly having open ends, that is, half-edges with only one vertex and a
boundary end-point. There is a simple operation, so-called {\it doubling},
that takes a one-trivalent graph $\Ga$ to a closed trivalent graph. It
consists simply of gluing two mirror copies of $\Ga$ along edges with
one-vertices:
 \begin{equation}
 \Ga \# \ov{\Ga}.
 \end{equation}
For example, the simplest graph of genus 2 without loops in the shape of 
$\Theta$ is the
double of the tristar, the one-trivalent graph having only one vertex.

This allows us to prescribe the value $V(\Ga_{j})$ of any colored
trivalent
graph $\Ga_j$, even if the function $j$ doesn't satisfy the
Clebsch--Gordan
conditions: removing all vertices not satisfying (1.6), we get a disjoint
union of one-trivalent connected components, that we can double to get a
disjoint union of spin networks. The sum of their values defines the value
of any colored trivalent graph $\Ga_j$.

The rest of the section comments briefly on the notion of spin network for
mathematicians.

\subsection{Why such an object?}
A spin network realizes a simple model of quantum geomery that is at the
same time discrete and purely combinatorial, and does not refer to any
background notion of space, time or geometry. The system consists of a
number of ``units'', each having a total angular momentum (a
representation
$j\in \widehat{SU(2)}$
 if the system has symmetry group $\SU(2)$). They interact in ways
that conserve the symmetry (or the total angular momentum). Thus any
interval
$e$ (edge) is a propagator of this unit; an event is an end (vertex) of
$e$
at which it meets two other edges whose labels satisfy (1.6). An
intertwiner
is necessarily nontrivial. A spin network is obtained by continuing in
this
way until we get a closed graph.

Hence such a model is described by an arbitrary trivalent graph with edges
labeled by integers (colors = twice the total angular momentum). The
vertices
describe the interactions. An orientation is just an orientation of the
procedure from past to future. A change of orientation is just a change of
a ``time direction''. A network with open ends (that is one-trivalent
graph)
is a history. The connected sum with its mirror image, that is, its double
is a spin network (for details and other physical ideas behind this notion
see \cite{7}).

\subsection{Why trivalent?}
Multivalent graphs can be reduced to the trivalent case as follows. Take
each
$n$-valent vertex and replace it by a $n$-leaved tree with trivalent
vertices.
This tree has $n-3$ new ``internal edges''. A basis of intertwining
operators
for the original vertex is then given by all labelings of these internal
edges by spins satisfying the Clebsch--Gordan condition. There are many
different trees with $n$ leaves, and thus many such bases. To change from
one
basis to another requires repeated use of some standard matrices from
recoupling theory based on the Racah sum rules and the Biedenharn--Elliott
identity in term of the $6j$ symbols (see \cite{8} and \cite{9}).

There are a number of {\it applications} of spin networks to
``continuous''
theories.

\subsection{Lattice gauge theory}
Spin networks are a generalization of knots and links if we consider them
as graphs embedded in space. They can be used in place of the regular
cubic
lattice in lattice models of gauge theories. Moreover, the central
kinematic
concept in quantum gravity is that the space of diffeomorphism-invariant
states is spanned by a basis in one-to-one correspondence with (orbits of) 
embeddings of
spin networks (see \cite{7}).

\subsection{More algebra}
The physical origin of spin networks dictates the labeling by
representations
of Lie groups. But in this definition, the main property that we need to
switch on is the following:
 \begin{quote}
 any {\it product} of two labels can be decomposed into a {\it sum} of
labels.
 \end{quote}
More precisely, the product of two labels defines a finite set of labels,
and
the Clebsch--Gordon condition is just a choice of an element of this set.
This algebraic structure is usually called a {\it category with tensor
product}.

There are algebras whose representation theory has this property; these
more general objects are Hopf algebras. The algebraic structure of
their representation theory can be described in terms of monoidal
categories. There is a still more general class of spin network associated
with these objects. It is also traditional to use representations of
quantum groups (deformations of Lie algebras) as labels. This is
technically
quite reasonable: they satisfy a modified set of recoupling identities --
quantum $6j$ symbols that depend on the parameter $q$ that can be
specified.
New invariants of 3-manifolds can be constructed using these \cite{10}.
However the labeling by representations of quantum groups differs from
ordinary spin networks in several ways. They do not correspond to groups,
and thus do not correspond to gauge invariants of classical connections.
(However the limit $q\to 1$ corresponds to the classical limit
$k\to\infty$
of Chern--Simons theory, sending the Kauffman bracket $K^{k}(\Ga)$ to the
Penrose value $V(\Ga)$. This reflects a deep mathematical relationship
between the representation theory of quantum groups at roots of unity and
the representation theory of the corresponding loop group of level $k$.)

Other geometric objects that could be used as labels are exceptional
bundles
\cite{11}. The many mathematical possibilities for labeling spin networks
stimulates approaches to realize the main expectation of experts in
physics:
 \begin{quote}
 there is a quantum theory ``$X$'', {\it defined purely algebraically},
not
involving any background  geometry, whose classical limit is
$3+1$ general relativity coupled to certain matter fields.
 \end{quote}

The theory ``$X$''  realizes directly the holographic
conjecture and the Bekenstein bound. We recall that 'tHooft and Susskind's
holographic conjecture states that such a theory is defined in terms of
state spaces and observables on surfaces (see for example \cite{12}). It
is
mathematically quite reasonable to develop a correspondence sending spin
networks to the geometry (topology) of surfaces; but first we describe the
second (well known) face of spin networks.

\section{Second face}
\subsection{Harmonic analysis}
The main motivation behind spin networks is to quantize general
relativity.
For this, we need a Hilbert space of states with a collection of operators
as observables. It is reasonable to expect that the states are functions
on
the configuration space of a system (or more generally, sections of a
bundle).
Following this thread, we start by sending our spin networks to functions
on
some space corresponding to a spin network, and prove later (in Section~3)
that the space is actually independent of the spin network. As before, we
construct this space using only the first component -- the trivalent
graph.
Here we need to use harmonic analysis on groups.

Consider the product
 \begin{equation}
 \SU(2)^{E(\Ga)}=\prod_{e\in E(\Ga)} \SU(2)_e
 \end{equation}
with $\SU(2)$ components enumerated by edges of $\Ga$, and the product
 \begin{equation}
 \SU(2)^{V(\Ga)}=\prod_{v\in V(\Ga)} \SU(2)_v
 \end{equation}
with components enumerated by the vertices. Let $dx$ be the Haar measure
on
$\SU(2)$ normalized by the condition $\int_{\SU(2)}dx=1$ and $\vec dx$ the
product measure on $\SU(2)^{E(\Ga)}$ normalized by $\int \vec dx=1$. Then
by
the Peter--Weyl formula, any function $f\in L^{2}(\SU(2)^{E(\Ga)}, \vec
dx)$
has the decomposition
 \begin{equation}
 f(x)=\sum_{\vec{\rho}\in \widehat{\SU(2)^{E(\Ga)}}}
 \Tr[B_{\vec{\rho},f} \vec{\rho}(x)],
 \end{equation}
where $\widehat{\SU(2)^{E(\Ga)}}$ is the space of irreducible
representations
of $\SU(2)^{E(\Ga)}$, and $B_{\vec{\rho},f}$ are endomorphisms of the
space
$V_{\vec{\rho}}$ of the representation $\vec{\rho}$, given by
 \begin{equation}
 B_{\vec{\rho},f}=\frac{1}{\dim V_{\vec{\rho}}}\int_{\SU(2)^{E(\Ga)}}
 f(x) \vec{\rho}\,^{-1}(x) \vec dx.
 \end{equation}
Recall that every irreducible representation of $\SU(2)^{E(\Ga)}$ is given
by tensor product of irreducible representations of $\SU(2)$:
 \begin{equation}
 \vec\rho=\rho_{1}\tensor \cdots \tensor\rho_{3g-3}.
 \end{equation}
This is of course an analog of the standard Fourier decomposition. Here a
representation $\vec\rho\,$ is a {\it label of a frequency} and an
endomorphism
$B_{\vec\rho, f}$ is the {\it Fourier coefficient}, that is, a {\it
number}.
The last formula is nothing other than the integral formula for a Fourier
coefficient.

Therefore every spin network $\Ga_{j}$ of genus $g$ defines a
representation
of $\SU(2)^{E(\Ga)}$ by the tensor product of all labels
 \begin{equation}
 \vec{\dotlessj}=\bigotimes_{e\in E(\Ga)} j_e
 \end{equation}
(the label of a frequency) and to get a function on $\SU(2)^{E(\Ga)}$ we
must
define an endomorphism $B(\Ga_j)$ using a labeling of a spin
network. However,
any endomorphism of the space $V_{\vec{\dotlessj}}$ is a vector in the
tensor product 
 \begin{equation}
 \Bigl(\bigotimes_{e\in E(\Ga)} j_e\Bigr)
 \tensor \Bigl(\bigotimes_{e\in E(\Ga)} j_e\Bigr)^*=
 \Bigl(\bigotimes_{e\in E(\Ga)} j_e\Bigr)
 \tensor \Bigl(\bigotimes_{e\in E(\Ga)} j_e\Bigr),
 \end{equation}
since we are dealing with $\SU(2)$-representations. But components of the
final product can be labeled by elements of the set $F(\Ga)$, and by
(1.2--3)
we can decompose it as
 \begin{equation}
 \Bigl(\bigotimes_{e\in E(\Ga)} j_e\Bigr)
 \tensor \Bigl(\bigotimes_{e\in E(\Ga)} j_e\Bigr)^*=
 \bigotimes_{v\in V(\Ga)}(j_{v,1} \tensor j_{v,2} \tensor j_{v,3}).
 \end{equation}
For every triple representations around a vertex $v\in V(\Ga)$ we have a
vector
 \begin{equation}
 i_v\in j_{v,1} \tensor j_{v,2} \tensor j_{v,3}
 \end{equation}
as in (1.7) and their tensor product gives us the vector
 \begin{equation}
 B(\Ga_j)=\bigotimes_{v\in V(\Ga)} i_v\in \End V_{\vec{\dotlessj}}\,.
 \end{equation}

As we saw, this endomorphism is an analog of a number -- a Fourier
coefficient.
But in some sense this number is an {\it integer}. Indeed, to construct
the
endomorphism, we use integer blocks of representations and its matrix must
be
an integer with respect to the multiplicative components of the
representation. 

\subsection{Fourier term with integer coefficient}
Hence a spin network defines a Fourier term with integer
coefficient. Moreover,
we may {\it identify} a spin network with this Fourier term with integer
coefficient, and vice versa. Of course, every spin network $\Ga_{j}$ as a
Fourier term with integer coefficient defines a state, that is, a function
 \begin{equation}
 f_{\Ga_{j}}(x)=\Tr[B(\Ga) \vec{\dotlessj}(x)]\in L^{2}(\SU(2)^{E(\Ga)}, d
\vec x).
 \end{equation}

To switch on the action of $\SU(2)^{V(\Ga)}$ (2.2) on the space
$\SU(2)^{E(\Ga)}$ (2.1), consider an orientation of $\Ga$, that is,
orientations of the edges such that every vertex has edges both ``in'' and
``out''. Such an orientation always exists. Now for every oriented edge
$e$
with ends $v_{\mathrm{in}}$ and $v_{\mathrm{out}}$, we set
 \begin{equation}
 g(t_{e})=g_{v_{\mathrm{in}}} \circ t_{e} \circ g_{\mathrm{out}}^{-1},
 \end{equation}
where
 \begin{align}
 g &=(g_{v_{1}}, \dots, g_{v_{2g-2}})\in \SU(2)^{2g-2}_{V(\Ga)}, \\
 t &=(t_{1}, \dots, t_{3g-3})\in \SU(2)^{3g-3}_{E(\Ga)}.
 \end{align}

As a standard result of harmonic analysis on groups, we get the following
\begin{prop}
 If the endomorphism $B(\Ga_j)$ (2.10) is intertwining, then the function
$f_{\Ga_{j}}$ (2.11) is invariant under the action (2.12):
 \begin{equation}
 f_{\Ga_{j}}\in L^2(\SU(2)^{E(\Ga)})^{\SU(2)^{V(\Ga)}}.
 \end{equation}
That is, $f_{\Ga_{j}}$ is a function on the homogeneous space
 \begin{equation}
 Q_{\Ga}=\SU(2)^{E(\Ga)} / \SU(2)^{V(\Ga)}.
 \end{equation}
\end{prop}

 Now $Q_{\Ga}$ does not depend on the choice of the orientation and on
the labeling $j$ of the graph, but the function $f_{\Ga_{j}}$ is {\it
equivalent\/} to the full spin network: the labeling is recognized by its
Fourier decomposition.

Thus the second face of every spin network $\Ga_{j}$ is an {\it
integer}. How
to add two such ``numbers''? How to multiply them? These operations are
related to the ``interaction'' of the combinatorial string ``events'' (see
below). 

The practical experiences behind such number theoretic intuition of spin
networks led to the technique of Feynman diagrams in perturbative quantum
field theories and quantum computers.

Many years ago Jacobi, with his theory of theta functions of one variable,
and
Riemann related number theory (a theory with {\it discrete objects}) with
complex analysis of one variable (a theory with {\it nondiscrete objects},
and
admitting many limits). What we want to do in the full program \cite{5} is
similar: to send the theory of spin networks as a combinatorial theory
 to
``continuous'' theories (the Chern--Simons and WZW theories). The first
step
in this is the {\it third face} of spin networks.

\section{Third face}
\subsection{From spin networks to surfaces}
There are two ways of transforming a spin network $\Ga_{j}$ embedded in
some
space $Y$ into a surface. The first is well known from the point of view
of
deformation quantization (see for example \cite{13}); let us call it the
{\it ribbon method}. We describe it briefly because \cite{13} does not
treat
it in our way. Our spin network is a generalization of a knot or link, and
it can be framed in the same vein. Any framed trivalent graph can be
lifted
to a ribbon by the same trick as a knot. Now our ribbon is an oriented
Riemann
surface $S$ with finite set of holes. We can apply the technique that is
well
known in the theory of framed graphs to get some topological invariants of
the
pair $S\subset Y$. Following this standard method of knot theory, we can
deform the representation theory of the algebra $\fsu(2)$ to the
representation theory of the quantum group $\fsu_{q}(2)$ and so on. Using
our
analogy between spin networks and integers, we can compare the theory of
spin
networks with quantum groups with the theory over function fields in
algebraic geometry. Everybody knows that this theory is easier and
simpler.

\subsection{Pumping up trick}
The ``compact'' method is much more interesting for us: we pump up the
edges
of $\Ga$ to tubes and the vertices to small 2-spheres. We get a Riemann
surface $\Si_{\Ga}$ of genus $g$ marked by a tube $\{\tilde{e}\}$ for
every
$e\in E(\Ga)$ and a {\it trinion} $\{ \tilde{v} \}$ for every $v\in
V(\Ga)$,
where each trinion is a 2-sphere with three holes. The isotopy classes of
meridian circles of tubes define $3g-3$ disjoint, non\-contractible,
pairwise nonisotropic classes. Let us consider any representations of
these
classes as simple loops (or circles) $\{C_{i}\}$ on $\Si$. The complement
is
the union
 \begin{equation}
 \Si_{\Ga} \setminus \{C_{1},\dots,C_{3g-3}\}=
 \bigcup_{i=1}^{2g-2} \tilde{v}_{i}
 \end{equation}
of $2g-2$ trinions corresponding to vertices of our graph $\Ga$.

Moreover, we can construct a map
 \begin{equation}
 m \colon \Si_{\Ga}\to \Ga
 \end{equation}
such that for every point $p$ in the open edge $m^{-1}(p)=S^{1}$, and for
every vertex $v$ the preimage $m^{-1}(v)=\text{``$\infty$''}$ (that is, a
bouquet of 2 circles). Thus the result of this pumping up can be viewed as
a {\it dynamics of circles} in $Y$.

This is the geometric idea of string dynamics, rather than the particle
dynamics described in 1.2. We call it {\it combinatorial string dynamics}.
In this natural picture an event or an interaction of two surfaces (or
string
families) is its intersection. A result of an interaction of two surfaces
$\Si', \Si''$ at a point $p\in \Si', \Si''$ is just the connected sum
 \begin{equation}
 \Si=\Si' \#_p \Si'',
 \end{equation}
and to get a trinion decomposition of $S$ we have to solve a local
problem.
Indeed, it is enough to consider the case when a point $p$ is contained in
two tubes $p\in \tilde{e'}\subset \Si', \quad p\in \tilde{e''}\subset
\Si''$.
Thus we only need to decompose a connected sum $\tilde{e'}\#_p
\tilde{e''}$.
The connected sum operation gives us a new tube $\tilde{p}$ (the ``neck''
of
the connected sum) with its meredian and the trinion decomposition of
$\tilde{e'}\#_p \tilde{e''}$. This procedure is well known and parallel to
the construction of a trivalent graph $\Ga$ from two trivalent graphs
$\Ga'$,
$\Ga''$ by joining the midpoints of two edges $m'\in e'\in \Ga'$ and
$m''\in e''\in \Ga''$ by a new edge $p$. We get two new trivalent vertices
$m'$ and $m''$ with edges $e'_1,e'_2,p$ around $m'$ and $e''_1,e''_2,p$
around $m''$. Pairs of edges $e'_1,e'_2$ and $e''_1,e''_2$ are ``halves''
of the edges $e'$ and $e''$. They inherit colors from these edges:
$j (e'_i)=j (e')$ and $j (e''_i)=j (e'')$. The final problem is to
prescribe the color of the last tube $\tilde{p}$. We can do this by
prescribing
any number $j(p)$ not contradicting the Clebsch--Gordan conditions around
$m'$ and $m''$. 

There are interpretations in this vein of the area operator (see \cite{2},
Fig.~2), the intersection of spin networks with a boundary surface in a
3-dimensional cylinder and so on, although we don't have space to go into
this here. (It is interesting to try to extract these combinatorial
objects
from nonperturbative string theories, for example from \cite{14}). 
 
This approach to constructing a good theory leads to an amazing geo\-metry
of
embedded graphs: a pumped up surface can itself be knotted, and we have to
develop a theory of 2-dimensional knots. Secondly there are others
versions 
writing down the amplitude of an event, that is, some intersection index
of
surfaces.

In the modern (super) string theorie \cite{15}, two surfaces interact via
a
metric on $Y$. In some dimensions, such theories do not have anomalies,
but
of course we have lost background-independence at the same time. The
relation
between this string theory and the expected background-independent string
theory is reflected in the quantum theory ``$X$'' at the end of
1.4. Namely the
perturbative theory around the classical limit of ``$X$''' must be
described
by the modern perturbative string theory.

\subsection{Representations space}
We stress again that up to now we have only used the trivalent graph
$\Ga$.
We forget for a moment the function $j \colon E(\Ga)\to \widehat{SU(2)}$
(1.5). Now
consider the space $R_{g}$ of gauge classes of flat $\SU(2)$-connections
on
$\Si$, that is, the space
 \begin{equation}
 R_{g}=\Hom(\pi_{1}(\Si_{\Ga}),\SU(2))/\PU(2)
 \end{equation}
of conjugacy classes of $\SU(2)$ representations of the fundamental
group of our Riemann surface. Our trinion decomposition (3.1) defines a
map
 \begin{equation}
 \pi_{\Ga} \colon R_{g}\to \R^{3g-3},
 \end{equation}
where the last space is just Euclidean space with the special coordinate
system
$(c_{1}, \dots,c_{3g-3})$. For a class of representations $\rho\in R_{g}$
 \begin{equation}
 c_{i}(\rho)=\frac{1}{\pi} \cdot \cos^{-1}\Bigl(\frac{1}{2}\Tr
 \rho([C_{i}])\Bigr)\subset [0,1].
 \end{equation}
It is well known that the functions $c_{i}$ on $R_{g}$ are continuous on
all
$R_{g}$ and smooth over $(0,1)$. Moreover, the image of $R_{g}$ under
$\pi_{\Ga}$ is a convex poly\-hedron
 \begin{equation}
 \De_{\Ga}\subset [0,1]^{3g-3}.
 \end{equation}
Now we use the second component $j$ of our spin network $\Ga_{j}$. For
this,
suppose that our spin network is of level $k$. Then for every label, we
can
consider the number $\frac{j_{e_{i}}}{k}$ as the $i$th coordinate of a
point
in $\R^{3g-3}$. So the function $\frac{j}{k}$ as coordinates defines a
point
$p_{j}\in \R^{3g-3}$. It was proved in \cite{16} that
 \begin{equation}
 p_{j}\in \De_{\Ga}.
 \end{equation}
So we get the subcycle
 \begin{equation}
 \pi_{\Ga}^{-1}(p_{j})\subset R_{g}
 \end{equation}
This is a serious geometric object deserving careful study.

First of all, the space $R_{g}$ is equipped with the canonical symplectic
form $\Om_{G}$. We call it the {\it Goldman form} \cite{17}. Thus the pair
$(R_{g}, \Om_{G})$ is a phase space of a classical mechanical system. This
system is completely integrable, and the complete set of first integrals
is
given by the coordinates $c_{i}$ in $\R^{3g-3}$. So the map $\pi_{\Ga}$
(3.5)
is a real polarization of this mechanical system, and the general fiber of
this map is a $(3g-3)$-dimensional Lagrangian torus.

Moreover, our phase space has canonical prequantization data
$(\Theta,a_{CS})$, where $\Theta$ is Hermitian line bundle with unitary
connection $a_{CS}$ (the so-called Chern--Simons connection), whose
curvature form satisfies
 \begin{equation}
 F_{a_{CS}}=2 \pi i \Om_{G}.
 \end{equation}

Recall that a Lagrangian cycle $S\subset R_{g}$ is a Bohr--Sommerfeld
cycle
of level $k$ iff the restriction $(\Theta^k, k a_{CS}) \vert_{S}$ admits a
covariant constant section (see \cite{18}). Let
 \begin{equation}
 \SNW_{k}(\Ga)=\{\Ga_{j}\}
 \end{equation}
be the set of all spin networks of level $k$ over a graph $\Ga$. Then we
have a set of Lagrangian fibers of the fibration (3.5):
 \begin{equation}
 \BS_{k}(\Ga)=\{\pi_{\Ga}^{-1}(p_{j})\}.
 \end{equation}

The following result was proved in \cite{16}:
\begin{prop} The collection $\BS_{k}(\Ga)=\{\pi_{\Ga}^{-1}(p_{j})\}$ is
the
set of all Bohr--Sommerfeld fibers of level $k$.
\end{prop}

Now we can fix the third face of a spin network $\Ga_{j}$ of genus $g$ and
level $k$, its geometric equivalent: this is a {\it Bohr--Sommerfeld fiber
of the real polarization of} $R_{g}$ {\it given by the trinion
decomposition
of Riemann surface} $\Si_\Ga$. Thus the number $\vert \BS_{k}(\Ga) \vert$
of
Bohr--Sommerfeld fibers of level $k$ is equal to the number
$ N^{k}(\Ga) $ of spin networks of level $k$ over a graph $\Ga$.
But this number is equal to the Verlinde number (see \cite{16} and
\cite{18}):
 \[
  N^{k}(\Ga) =\vert \BS_{k}(\Ga) \vert=\frac{(k+2)^{g-1}}{2^{g-1}} 
\sum_{n=1}^{k+1}
 \frac{1}{(\sin(\frac{n\pi}{k+2}))^{2g-2}}\,.
 \]

\subsection{More ``bulk'' geometry}

In the same vein as $\Si_{\Ga}$, our trivalent graph $\Ga$ defines a
{\it handlebody } $\widetilde{\Si_{\Ga}}$, that is, a 3-manifold with
boundary $\Si_{\Ga}$. We have the epimomorphism
 \begin{equation}
 r \colon \pi_1(\Si_{\Ga})\to \pi_1(\widetilde{\Si_{\Ga}}).
 \end{equation}
The fundamental group $\pi_1(\widetilde{\Si_{\Ga}})$ is the free group
with
generators, say $\{b_i\}$ for $i=1,\dots,g$. We can consider them as
elements
of $\pi_1(\Si_{\Ga})$. And the kernel of $r$ is the free group with
generators, say $\{a_i\}$ for $i=1,\dots,g$. Then
 \begin{equation}
 \pi_1(\Si_\Ga)=\Span{a_1,\dots, a_g, b_1,\dots, b_g \Biggm|
 \prod_{i=1}^{g} [a_i, b_i]=\id}.
 \end{equation}
is the standard presentation of the fundamental group.

Starting with a trivalent graph $\Ga$ we get a Riemann surface $\Si_{\Ga}$
with the collection of elements $\{[C_{i}]\}$ of $\pi_1(\Si_\Ga)$
corresponding to disjoint loops (3.1). Between these elements we can fix
the
first $g$ loops to be $[C_{1}]=a_1, \dots, [C_{g}]=a_g$ by changing the
numbering, so that we can add elements $b_{1},\dots,b_{g}\in\{C_{i}\})$ to
make a standard basis (3.14) of the fundamental group $\pi_1(\Si_{\Ga})$.

A graph with such additional choice is a {\it marked} graph $\Ga^{m}$.
Returning to the space $R_{g}$ (3.4), we get the subcycle
 \begin{equation}
 uS_{g}=\{\rho\in R_{g} \bigm| \rho([C_{i}])=1 \text{ for } i=1,\dots,g\},
 \end{equation}
that we call the {\it unitary Schottky space} of genus $g$.

This space can be presented as a homogeneous space
 \begin{equation}
 uS_{g}=\SU(2)^{g} / \Ad_{\mathrm{diag}} \SU(2)
 \end{equation}
where $\Ad_{\mathrm{diag}}\SU(2)$ is the diagonal adjoint action on the
direct product.

The following statements are proved in \cite{16} (see also \cite{18})

\begin{prop} The unitary Schottky space (3.15) is a fiber of the real
polarization $\pi_{\Ga^{m}}$ (3.5). More precisely
 \begin{equation}
 uS_{g}=\pi_{\Ga^{m}}^{-1}(0,\dots, 0).
 \end{equation}
In particular
 \begin{equation}
 uS_{g}\in \BS_{k}(\Ga)=\SNW_{k}(\Ga)
 \end{equation}
for every level $k$. (See (3.11) and (3.12)).
\end{prop}

The interpretation of the space $R_{g}$ as the space of gauge classes of
flat connections on $\Si_{\Ga}$ gives the description of all fibers of
$\BS_{k}(\Ga)$: let $\pi_{\Ga}^{-1}(p_{j})$ be such a fiber. For each loop
$C_{i}$ for $i=1,\dots,3g-3$, write
 \begin{equation}
 Z_{e_{i}}^{j}=\text{stabilizer of } \pi_{\Ga}^{-1}(p_{j}) \vert_{C_{i}}
 \end{equation}
and
\begin{equation}
Z_{v_{i}}^{j}=\text{stabilizer of } \pi_{\Ga}^{-1}(p_{j})
\vert_{P_{i}}\end{equation}
where $P_{i}$ is a trinion of the decomposition (3.1).

We have two groups
 \begin{align}
 Z_{E(\Ga)}^{j}&=\prod_{e\in E(\Ga)} Z_{e}^{j} \\
 Z_{V(\Ga)}^{j}&=\prod_{v\in V(\Ga)} Z_{v}^{j}
 \end{align}
and the action of the second group on the first given by the formula
(2.12). Then the fiber is
 \begin{equation}
 \pi_{\Ga}^{-1}(p_{j})=Z_{E(\Ga)}^{j} / Z_{V(\Ga)}^{j}.
 \end{equation}
For general $j$
 \begin{align}
 Z_{e}^{j}&=\U(1) \text{ for every } e\in E(\Ga) \\
 Z_{v}^{j}&=Z_{2} \text{ is the center of } \SU(2)
 \end{align}
Thus for general $p_{j}\in \BS_{k}(\Ga)$ the fiber
 \begin{equation}
 \pi_{\Ga}^{-1}(p_{j})=\U(1)^{3g-3}
 \end{equation}
is a $(3g-3)$-torus.

 But for the special point,
 \begin{equation}
 \pi_{\Ga}^{-1}(0, \dots,0)=\SU(2)^{3g-3}_{\Ga} /
 \SU(2)^{2g-2}_{\Ga}=Q_{\Ga}
 \end{equation}
(see (3.23)), since
 \begin{equation}
 \pi_{\Ga}^{-1}(0,\dots,0) \vert_{C_{i}}=1 \quad\text{for } i=1,\dots,3g-3
 \end{equation}
and
 \begin{equation}
 \pi_{\Ga}^{-1}(0,\dots,0)\vert_{P_{i}}=1, \quad\text{for }
i=1,\dots,2g-2.
\end{equation}

Thus we have the identification
 \begin{equation}
 uS_{g}=Q_{\Ga}.
 \end{equation}
Of course it depends on the marking (3.15) of $\Ga$. Comparing the
probabilistic measures on components of products (2.1), (2.2) and
(3.16) gives
the following result:

\begin{prop} Under the identification (3.30)
 \begin{equation}
 f_{\Ga_{j}}\in L^{2}(\SU(2)^{g}, d \vec x)^{\Ad_{\mathrm{diag}}\SU(2)}
 \end{equation}
(see (2.15)).
\end{prop}

Thus we have collected all spin network states as functions on the same
space, that is, on the unitary Schottky space
 \begin{equation}
 f_{\Ga_{j}}\in L^{2}(sU_{g}, \vec dx).
 \end{equation}
This is just what we need to come to the general theory \cite{19}. Here we
send $g$ to infinity to get all differential invariant states. On the
other
hand, we have a bridge to well known theories \cite{20} and \cite{21}. The
following partial case of the constructions under consideration can be
added
to \cite{22}.

\section{Illustration: the Abelian case}

For $\U(1)$-spin networks 
 \begin{equation}
 \widehat{U(1)}=\Z\hidot,
 \end{equation}
and the triangle inequality (1.6) becomes the equality
 \begin{equation}
 j_{v,1}+j_{v,2}+j_{v,3}=0.
 \end{equation}

The harmonic analysis is just the classical Fourier decomposition. Now
 \begin{equation}
 \U(1)^{E(\Ga)}=\prod_{e\in E(\Ga)} \U(1)_e
 \end{equation}
and 
 \begin{equation}
 \U(1)^{V(\Ga)}=\prod_{v\in V(\Ga)} \U(1)_v
 \end{equation}

It is easy to see that under the action (2.12), the diagonal of
$\U(1)^{V(\Ga)}$ acts trivially. Thus in this case
 \begin{equation}
 Q_{\Ga}=\U(1)^{E(\Ga)} / \U(1)^{V(\Ga)}=\U(1)^g.
 \end{equation}

Now the representation space (3.4)
 \begin{equation} J_{\Si_\Ga}=
 \Hom(\pi_{1}(\Si_{\Ga}), \U(1)) 
 \end{equation}
is the Jacobian of our surface. This group is the direct product of $2g$
copies of $\U(1)$, but components can be labeled by the basis (3.14)
 \begin{equation}
 J_{\Si_\Ga}=\prod_{i=1}^g \U(1)_{a_i} \times \prod_{i=1}^g \U(1)_{b_i}
 \end{equation}

We can view the coordinates (3.6) of the map (3.5) just as elements of 
target group $\U(1)$. Thus we have the map (3.5) in this case 
 \begin{equation}
 \pi_{\Ga} \colon J_{\Si_\Ga}\to \U(1)^{E(\Ga)}.
 \end{equation}

It is easy to see that the image is a $g$-torus
 \begin{equation}
 \De_{\Ga}=T^g_-\subset \U(1)^{E(\Ga)}
 \end{equation}
such that the projection
 \begin{equation}
 \pr \colon \U(1)^{E(\Ga)}\to \prod_{i=1}^g \U(1)_{c_i} 
\end{equation}
defines the isomorphism of $T^g_-$ to $\prod_{i=1}^g \U(1)_{c_i}$. 

Recall that we have the identification $a_i=c_i$ (3.14). Thus the map
$\pi_\Ga$ (3.5) is just the projection of the direct product (4.7) to
the component
 \begin{equation}
 \pi_\Ga\colon\prod_{i=1}^g\U(1)_{a_i}\times\prod_{i=1}^g\U(1)_{b_i}\to
 \prod_{i=1}^g\U(1)_{a_i}=T^g_-.
 \end{equation}

But the intersections of 1-cycles on our Riemann surface defines the
integral
symplectic form $\Om$ on $J_{\Si_\Ga}$, that is, the polarisation of the
Jacobian. It is easy to see that the fibration (4.11) is Lagrangian, that
is,
each fiber is a Lagrangian torus. 

Now consider a $\U(1)$-spin network of level $k$, that is, for every $e\in
E(\Ga)$, $0 \leq j_e \leq k-1$. Then the function $\frac{j}{k}$ defines a
point $p_j\in \U(1)^{E(\Ga)}$ and moreover $p_j\in \De_{\Ga}=T^g_-$. This
point is a point of order $k$ on the torus $T^g_-$. 

It is easy to see that the function $j$ (1.5) can be reconstructed from
this
point and every point of order $k$ of $T^g_-$ is an image of some
commutative
spin network. To say nothing of the fact that all fibers over these points
are Bohr--Sommerfeld fibers of the real polarization $\pi_\Ga$ (see
\cite{18}).
 
Thus the number $\vert \BS_{k}(\Ga) \vert$ of Bohr--Sommerfeld fibers of
level $k$ is equal to number $ N^{k}_A(\Ga) $ of Abelian spin networks
of level $k$ over a graph $\Ga$ and is equal to the number of points of
order
$k$ on $g$-torus $T^g_-$. But this number is equal to (see  \cite{18}):
 \[
  N^{k}_A(\Ga)=\vert\BS_{k}(\Ga)\vert=\vert(T^g_-)_k\vert=k^g.
 \]
The full program (all the stones of the mosaic) of the theory of classical
theta functions for Abelian case is developed in \cite{23}. 

\subsection*{Acknowledgments}

I would like to express my gratitude to my collaborators C. Florentino,
J. Mour\~ao, J.P. Nunes and to the Instituto Superior Tecnico of Lisbon
for
support and hospitality. Special thanks to Miles Reid for his permanent
help.

\bigskip
\noindent
\small{Algebra Section, Steklov Math Institute, Ul.\ Gubkina \\
Moscow 117333, Russia \\
email: Tyurin@tyurin.mian.su \\
{\it or} Tyurin@Maths.Warwick.Ac.UK }

\end{document}